\documentclass[12]{article}
\usepackage{amsfonts}
\usepackage{amsmath}
\usepackage[unicode,colorlinks,plainpages=false]{hyperref}
\usepackage{theorem}
\newtheorem{theorem}{Theorem}

\newtheorem{example}[theorem]{Example}

\newcommand{\openbox}{$\begin{array}{c}
\hspace*{-0.55em}\sqcap \hspace*{-0.60em}\\[-0.4em] \hline
\multicolumn{1}{c}{\hspace*{-0.60em}}\\[-0.8em]
\end{array}
$}

\usepackage{enumerate}

\begin{document}

\centerline{{\bf REMARKS ON THE PAPER}\footnote{This version differs from the published one that we clarified some part of the proof of Theorem~\ref{main}, and corrected the typing errors. {\bf Key words}: semigroups, left reductive semigroups, congruences;
{\bf AMS Classification}: 20M10; {\bf e-mail}: nagyat@math.bme.hu}}
\centerline{\bf "M. KOLIBIAR, ON A CONSTRUCTION OF SEMIGROUPS"}
\bigskip
\centerline{Attila Nagy}
\medskip
\centerline{Department of Algebra}
\centerline{Budapest University of Technology and Economics}
\centerline{P.O. Box 91}
\centerline{1521 Budapest}
\centerline{Hungary}

\bigskip

\begin{abstract} In his paper "On a construction of semigroups", M. Kolibiar gives a construction for a semigroup $T$ (beginning from a semigroup $S$) which is said to be derived from the semigroup $S$ by a $\theta$-construction. He asserted that every semigroup $T$ can be derived from the factor semigroup $T/\theta (T)$ by a $\theta$-construction, where $\theta (T)$ is the congruence on $T$ defined by: $(a, b)\in \theta (T)$ if and only if $xa=xb$ for all $x\in T$. Unfortunately, the paper contains some incorrect part. In our present paper we give a revision of the paper.
\end{abstract}

\medskip

A semigroup $S$ is called a left reductive semigroup (\cite{Clifford1:sg-1}, \cite{Nagybook:sg-2}) if, for every $a, b\in S$, whenever $xa=xb$ holds for all $x\in S$ then $a=b$. It is known that, for an arbitrary semigroup $S$, the relation $\theta (S)$ defined by $(a, b)\in \theta (S)$ for some $a, b\in S$ if and only if $xa=xb$ for all $x\in S$ is a congruence on $S$.
In \cite{Kolibiar:sg-4}, the author examined this congruence. He defined a sequence $\theta _n$ ($n=1, 2, \dots $) of congruences on an arbitrary semigroup $S$ as follows: $\theta _1=\theta (S)$, and if $\theta _n$ is given, $\theta _{n+1}$ is the congruence relation on $S$, induced by the congruence relation $\theta (S)/\theta _n$ (\cite{Howie:sg-5}). In Lemma 2 of \cite{Kolibiar:sg-4} it is shown that $(a, b)\in \theta _n$ for some $a, b\in S$ if and only if $xa=xb$ for all $x\in S^n$.

In Theorem 1 of \cite{Kolibiar:sg-4}, it is asserted that $\theta ^*=\cup_{n=1}^{\infty}\theta _n$ is the least element in the set of all congruence relations $\theta$ on $S$ such that $S/\theta$ is a left reductive semigroup. The proof of this theorem is not correct. The author asserts that, from the result $(xa, xb)\in \theta ^*$ for all $x\in S$, it follows that $txa=txb$ for all $x\in S$ and all $t\in S^k$ for some $k\in N$ (that is
$(xa, xb)\in \theta _k$ for all $x\in S$). This is not correct, because $(xa, xb)\in \theta^*$ for all $x\in S$ means that, for every $x\in S$, there is a positive $n_x$ such that $(xa, xb)\in \theta _{n_x}$. From this it does not follow that there is a positive integer $k$ such that
$(xa, xb)\in \theta _k$ for all $x\in S$.

Theorem 2 of \cite{Kolibiar:sg-4} is also incorrect. In part $(a)$ of Theorem 2, the following construction is given: For a semigroup $S$ and each $x\in S$, let $T_x\neq \emptyset$ be a set such that $T_x\cap T_y=\emptyset$ if $x\neq y$. Let a mapping $f^y_x: T_x \mapsto T_y$ be given for all couples $x, y\in S$ such that $y=xu$ (in the paper $y=ux$, but it is a missprint) for some $u\in S$. Suppose further $f_y^z\circ f_x^y =f_x^z$. Given
$a, b\in \cup _{x\in S}T_x=T$, set $a\circ b=f_x^{xy}(a)$, where $a\in T_x$, $b\in T_y$. Then $(T; \circ )$ is a semigroup (this semigroup is said to be derived from $S$ by a $\theta$-construction), and each set $T_x$ is contained in a $\theta$-class of $T$. If $S$ is left reductive then the $\theta$-classes of $T$ are exactly the sets $T_x$ ($x\in S$). The assertions of this part of the theorem is correct. In part $(b)$ of Theorem 2, it is asserted that, for an arbitrary semigroup $T$, if $S$ denotes the factor semigroup $T/\theta (T)$ and $T_{[x]\theta}$ denotes the $\theta (T)$-class $[x]\theta$ of $T$ containing $x$, then $f_{[x]\theta}^{[xy]\theta}: a\mapsto ab$ ($a\in [x]\theta$, $b\in [y]\theta$) is a mapping on $[x]\theta$ to $[xy]\theta$ and, for all $c, d\in T$, $cd=f_{[c]\theta}^{[cd]\theta}(c)$. This last equation means that $cd=c\circ d$ for all $c, d\in S$. The next example shows that the assertion $cd=f_{[c]\theta}^{[cd]\theta}(c)$ for all $c, d\in T$ is not correct.

\medskip

\begin{example}\label{ex1}\rm
Let $T=\{e, a, u, v, 0\}$ be a semigroup defined by Table~\ref{Fig1}.

\begin{table}[htbp]
\begin{center}
\begin{tabular}{l|l l l l l }
$\cdot$&$e$&$a$&$u$&$v$&$0$\\ \hline
$e$&$e$&$a$&$0$&$0$&$0$\\
$a$&$a$&$e$&$0$&$0$&$0$\\
$u$&$u$&$v$&$0$&$0$&$0$\\
$v$&$v$&$u$&$0$&$0$&$0$\\
$0$&$0$&$0$&$0$&$0$&$0$\\
\end{tabular}
\caption{}\label{Fig1}
\end{center}
\end{table}

Let $\theta$ denote the congruence $\theta (T)$. It is easy to see that the factor semigroup $S=T/\theta $ is a left reductive semigroup. We note that $S$ is a semigroup which can be obtained from a two-element group $\{[e]\theta, [a]\theta\}$ by adjunction the zero $[0]\theta$.

Consider the semigroup $(T;\circ)$ which is derived from $T/\theta $ by a $\theta$-construction (see (b) of Theorem 2 of \cite{Kolibiar:sg-4}). As $[0]\theta =[0]\theta [0]\theta$, $[0]\theta =[0]\theta [e]\theta$, $[0]\theta =[0]\theta [a]\theta$, we have three possibility to define a mapping $f^{[0]\theta}_{[0]\theta}$ from $[0]\theta$ into $[0]\theta$.

In case $f_{[0]\theta }^{[0]\theta } \ x\mapsto x\cdot 0=0$ for all $x\in [0]\theta$, we have
\[u\circ e=f_{[u]\theta}^{[u\cdot e]\theta}(u)=f_{[0]\theta}^{[0]\theta}(u)=0\neq u=u\cdot e.\]

In case $f_{[0]\theta }^{[0]\theta } \ x\mapsto x\cdot e$ for all $x\in [0]\theta$, we get
\[u\circ a=f_{[u]\theta}^{[u\cdot a]\theta}(u)=f_{[0]\theta}^{[0]\theta}(u)=u\cdot e=u\neq v=u\cdot a.\]

In case $f_{[0]\theta }^{[0]\theta } \ x\mapsto x\cdot a$ for all $x\in [0]\theta$, we have
\[u\circ e=f_{[u]\theta}^{[u\cdot e]\theta}(u)=f_{[0]\theta}^{[0]\theta}(u)=u\cdot a=v\neq u=u\cdot e.\]

Under any possible choosing of the mapping $f_{[0]\theta }^{[0]\theta }$ we have a contradiction. Thus $T$ can not be derived from $T/\theta$ by a $\theta$-construction.
\end{example}
\bigskip

In the next we formulate a new and correct version of Theorem 2 of \cite{Kolibiar:sg-4}.

\begin{theorem}\label{main} (a) Let $S$ be a semigroup. For each $x\in S$, associate a set $T_x\neq \emptyset$. Assume that $T_x\cap T_y=\emptyset$ for all $x\neq y$. Assume that, for each triple $(x, y, xy)\in S\times S\times S$ is associated a mapping $f_{(x , y, xy)}$ of $T_x$ into $T_{xy}$ acting on the left.
Suppose further
$f_{(xy, z, xyz)}\circ f_{(x, y, xy)}=f_{(x, yz, xyz)}$ for all triples $(x, y, z)\in S\times S\times S$. Given $a, b\in \cup _{x\in S}T_x=T$, set $a\circ b=f_{(x, y, xy)}(a)$, where $a\in T_x$ and $b\in T_y$. Then $(T; \circ )$ is a semigroup, and each set $T_x$ is contained in a $\theta (T)$-class of $T$. If $S$ is left reductive then the $\theta (T)$-classes of the semigroup $T$ are exactly the sets $T_x$ ($x\in S$).

(b) Let $(T; \cdot)$ be a semigroup. Let $\theta$ denote the congruence $\theta (T)$, and let $S=T/\theta$. For an element $[x]\theta\in S$, let $T_{[x]\theta}$ be the $\theta$-class $[x]\theta$ of $T$. For arbitrary triple
$([x]\theta , [y]\theta, [x\cdot y]\theta )\in S\times S\times S$, let $f_{([x]\theta, [y]\theta, [x\cdot y]\theta)} :\ a\mapsto a\cdot b$ ($a\in [x]\theta$), where $b$ is an arbitrary element of $[y]\theta$. (We note that $a\cdot b=a\cdot [b]\theta$ for all $a, b\in T$. We also note that we consider all of the mappings
$f_{([x]\theta ,[z]\theta , [x\cdot z]\theta )}$, where $[z]\theta \in S$ satisfies $[x]\theta[y]\theta =[x]\theta[z]\theta$.)
For all $c, d\in T$, let $c\circ d=f_{([c]\theta, [d]\theta, [c\cdot d]\theta)}(c)$. Then $(T;\circ)$ is a semigroup which is isomorphic to the semigroup $(T;\cdot)$. Consequently, every semigroup $(T;\cdot)$ is isomorphic to a semigroup derived from the semigroup $S=T/\theta$ using the construction in part $(a)$ of the theorem.
\end{theorem}

\noindent
{\bf Proof}. (a): Let $(x, y, z)\in S\times S\times S$ be an arbitrary triple and let $a\in T_x$, $b\in T_y$, $c\in T_z$ be arbitrary elements. Then
\[a\circ (b\circ c)=\]
\[=a\circ f_{(y, z, yz)}(b)=f_{(x, yz, xyz)}(a)=(f_{(xy, z, xyz)}\circ f_{(x, y, xy)})(a)=(f_{(x, y, xy)}(a))\circ c=\]
\[=(a\circ b)\circ c.\]
Thus $(T; \circ )$ is a semigroup. Let $x\in S$ be arbitrary. We show that $T_x$ is in a $\theta(T)$-class of $T$. If $a, b\in T_x$ are arbitrary elements then, for every $y\in S$ and for all $t\in T_y$,
\[t\circ a=f_{(y, x, yx)}(t)=t\circ b\] which implies that $(a, b)\in \theta (T)$.

Assume that $S$ is left reductive. Let $T_x$ and $T_y$ be sets such that they are in the same $\theta (T)$-class of $T$. Let $a\in T_x$ and $b\in T_y$ be arbitrary elements. Then, for every $z\in S$ and $t\in T_z$, $t\circ a=t\circ b\in T_{zx}\cap T_{zy}$. Thus $zx=zy$ for every $z\in S$. As $S$ is left reductive, we get $x=y$. Hence the $\theta (T)$-classes of $T$ are exactly the sets $T_x$ ($x\in S$).

(b): It is clear that, for every $x, y\in T$, $f_{([x]\theta, [y]\theta, [x\cdot y]\theta)}$ maps $[x]\theta=T_{[x]\theta}$ into $T_{[x\cdot y]\theta}=[x\cdot y]\theta$.
For every $(x, y, z)\in T\times T\times T$ and $a\in [x]\theta$,
\[(f_{([x\cdot y]\theta , [z]\theta, [x\cdot y\cdot z]\theta )}\circ f_{([x]\theta, [y]\theta, [x\cdot y]\theta )})(a)=\]
\[=(a\cdot y)\cdot z=a\cdot (y\cdot z)=f_{([x]\theta, [y\cdot z]\theta, [x\cdot y\cdot z]\theta)}(a).\]
Thus $(T ;\circ)$ is a semigroup by part $(a)$ of the theorem.

For every $x,y\in T$ and every $a\in [x]\theta$, $b\in [y]\theta$,
\[a\circ b= f_{([x]\theta, [y]\theta, [x\cdot y]\theta)}(a)=a\cdot b.\] Thus $(T; \cdot )$ is isomophic to the semigroup $(T; \circ )$.\hfill\openbox



\begin{example}\label{ex2}\rm
Let $T$ be the semigroup defined in Example 1.  Let $\theta$ denote the congruence $\theta (T)$. Then $S=T/\theta$ is a semigroup which can be obtained from a two-element group by adjunction of a zero. Let $e, a, 0$ denote the $\theta$-classes $[e]\theta =\{ a\}$, $[a]\theta =\{ a\}$, $[0]\theta =\{ 0, u, v\}$ of $T$, respectively. Then the Cayley-table of $S$ is
\begin{table}[htbp]
\begin{center}
\begin{tabular}{l|l l l }
$\empty$&$e$&$a$&$0$\\ \hline
$e$&$e$&$a$&$0$\\
$a$&$a$&$e$&$0$\\
$0$&$0$&$0$&$0$\\
 \end{tabular}
\caption{}\label{Fig2}
\end{center}
\end{table}

Apply our construction when the beginning semigroup is $S$. The sets $T_x$ ($x\in S$) are $T_e=\{ e\}$, $T_a=\{ a\}$, $T_0=\{ 0, u, v\}$. The mappings (see $(b)$ of Theorem~\ref{main}) are:
$f_{(e, e, e)}:\ e\mapsto e$, $f_{(e, a, a)}:\ e\mapsto a$, $f_{(e, 0, 0)}:\ e\mapsto 0$, $f_{(a, e, a)}:\ a\mapsto a$, $f_{(a, a, e)}:\ a\mapsto e$,
$f_{(a, 0, 0)}:\ a\mapsto 0$, $f_{(0, e, 0)}:\ x\mapsto x\cdot e$ for every $x\in [0]\theta$, $f_{(0, a, 0)}:\ x\mapsto x\cdot a$ for every $x\in [0]\theta$,
$f_{(0, 0, 0)}:\ x\mapsto 0$ for every $x\in [0]\theta$. It is easy to see that the semigroup $T$ is isomorphic to the semigroup $(T; \circ )$ derived from $S$ by applying our construction in Theorem~\ref{main}.
\end{example}

\medskip

We formulate the new version of Theorem 3 of \cite{Kolibiar:sg-4}.

\begin{theorem}\label{sequence} Let $S$ be a semigroup whose congruence relations satisfy the ascending chain condition. Then there are semigroups $S_0=S, S_1, \dots S_n$ such that $S_n$ is left reductive, and each $S_{i-1}$ can be derived from $S_{i}$ ($i=1, \dots , n$) by the construction defined in $(a)$ of Theorem~\ref{main}.
\end{theorem}

\noindent
{\bf Proof}. Let $S$ be a semigroup. For a congruence $\varrho$ on $S$, consider the congruence $\varrho ^*$ defined by $(a, b)\in \varrho ^*$ for some $a, b\in S$ if and only if $(xa, xb)\in \varrho$ for all $x\in S$. (see \cite{Nagyleft:sg-3}). Let
\[\varrho ^{(0)}\subseteq \varrho ^{(1)}\subseteq \cdots \subseteq \varrho ^{(n)}\subseteq \dots\]
be a sequence of congruences on $S$ defined by $\varrho ^{(0)}=\varrho$ and, for every nonnegative integer $i$, let $\varrho ^{(i+1)}=(\varrho ^{(i)})^*$ (see \cite{Nagyleft:sg-3}). By Lemma 2 of \cite{Kolibiar:sg-4}, $\iota _S ^{(i)}=\theta _i$ for every positive integer $i$, where $\iota _S$ denotes the identity relation on $S$ and $\theta _i$ are the congruences examined above. As the congruence relations on $S$ satisfy the ascending chain condition, there is a (least) nonegative integer $n$ such that $\iota ^{(n)}_S=\iota ^{(n+1)}_S=\cdots$. Let $S_i=S/\iota ^{(i)}_S$ for every $i=0, 1, \dots , n$. Thus $S=S_0$.
By Theorem 1 of \cite{Nagyleft:sg-3}, $\iota ^{(n)}_S=\theta _n$ is a left reductive congruence which means that the factor semigroup $S_n=S/\theta _n$ is left reductive. For every $i=1, \dots , n$,
\[S_i=S/\iota ^{(i)}_S\cong (S/\iota ^{i-1}_S)/(\iota ^{(i)}_S/\iota ^{(i-1)}_S)\cong S_{i-1}/(\iota ^{(i)}_S/\iota ^{i-1)}_S)\cong S_{i-1}/((\iota ^{(i-1)}_S)^*/\iota ^{(i-1)}_S)\]
using also Theorem 5.6 of \cite{Howie:sg-5}.
By Lemma 7 of \cite{Nagyleft:sg-3},
\[(\iota ^{(i-1)}_S)^*/\iota ^{(i-1)}_S=(\iota_{S_{i-1}})^*=\theta (S_{i-1}).\] Thus \[S_i\cong S_{i-1}/\theta (S_{i-1})\] for every $i=1, \dots , n$. By $(b)$ of Theorem~\ref{main}, the semigroup $S_{i-1}$ ($i=1, \dots , n$) can be derived from $S_i$ by the construction in $(a)$ of Theorem~\ref{main}.\hfill\openbox

\bigskip

{\bf An addendum}: It may be that the semigroup $S_n=S/\theta _n$ in Theorem~\ref{sequence} has only one element, that is, $\theta _n$ is the universal relation $\omega _S$ on $S$. The next theorem is about this case.

\begin{theorem}\label{adendum} For a semigroup $S$, $\theta _n=\omega _S$ for some positive integer $n$ if and only if $S$ is an ideal extension of a left zero semigroup by a nilpotent semigroup.
\end{theorem}

\noindent
{\bf Proof}. Let $\iota _S$ denote the identity relation on a semigroup $S$. By Theorem 5 of \cite{Nagyleft:sg-3}, $\iota ^{(n)}_S=\omega _S$ for some non-negative integer $n$ if and only if $S$ is an ideal extension of a left zero semigroup by a nilpotent semigroup. As $\theta (S)=\iota ^{(1)}_S$,
our assertion is a consequence of Theorem 5 of \cite{Nagyleft:sg-3}.\hfill\openbox


\bigskip

\begin{thebibliography}{1}

\bibitem{Clifford1:sg-1} A.H. Clifford and G.B. Preston, {\it The Algebraic Theory of Semigroups I}, American Mathematical Society, Providence, R. I., 1961

\bibitem{Howie:sg-5} J.M. Howie, {\it An Introduction to Semigroup Theory}, Academic Press, London, 1976

\bibitem{Kolibiar:sg-4} M. Kolibiar,  {\it On a construction of semigroups}, Scripta Fac. Sci. Nat. Ujep Brunensis, Arch. Math., 22:99-100, 1971

\bibitem{Nagybook:sg-2} A. Nagy, {\it Special Classes of Semigroups}, Kluwer Academic Publishers, Dordrecht, Boston, London, 2001

\bibitem{Nagyleft:sg-3} A. Nagy, {\it Left reductive congruences on semigroups}, Semigroup Forum, 87:129-148, 2013

\end{thebibliography}
\end{document}